\documentclass[11pt]{amsart}
\usepackage{amscd} 
\usepackage{amsfonts} 
\usepackage{amssymb} 
\usepackage{latexsym}

\newcommand{\ncm}{\newcommand}

\newtheorem{theorem}{Theorem}[section]
\newtheorem{prop}[theorem]{Proposition}
\newtheorem{lemma}[theorem]{Lemma}
\newtheorem{cor}[theorem]{Corollary}
\newtheorem{lem&def}[theorem]{Lemma \& Definition}
\newtheorem{definition}[theorem]{Definition}
\newtheorem{example}[theorem]{Example}

\def\C{\mathbb{C}\,} 
\def\Z{\mathbb{Z}\,} 
 
\def\H{\mathbb{H}\,}
\def\Q{\mathbb{Q}\,}

\def\id{\mbox{\rm id}}

\def\into{\hookrightarrow}
\def\to{\rightarrow}

\def\End{\mbox{\rm End}\,}

\def\Hom{\mbox{\rm Hom}\,}
\def\o{{\otimes}}    %tensor product 
\def\|{\, | \, }
\def\bra{\langle}
\def\ket{\rangle}

\ncm{\rarr}[1]{\stackrel{#1}{\longrightarrow}}
\ncm{\larr}[1]{\stackrel{#1}{\longleftarrow}}

\def\cop{\Delta}

\def\eps{\varepsilon}

\def\du1{\hat 1}

\def\-1{_{(-1)}}
\def\0{_{(0)}}
\def\1{_{(1)}}
\def\2{_{(2)}}
\def\3{_{(3)}}

\def\du1{\hat 1}
\def\lact{\triangleright}
\def\ract{\triangleleft}

\begin{document}

\title[Finite depth and a Galois correspondence]{Finite depth and Jacobson-Bourbaki correspondence}
\author{Lars Kadison} 
\address{Department of Mathematics \\ University of Pennsylvania \\
David Rittenhouse Lab, 209 S. 33rd St. \\ 
Philadelphia, PA 19104} 
\email{lkadison@math.upenn.edu}
%\URL{www.math.upenn.edu/~lkadison} 
% author two information 
\thanks{To David Harbater and the Penn Galois seminar  thanks for inviting me to give a talk and posing several  
interesting questions.}
\date{} 

\begin{abstract}
 We introduce a notion of 
depth three  tower of three rings $C \subseteq B \subseteq A$
with depth two ring extension $A \| B$ recovered when $B = C$. If $A = \End B_C$ and $B \| C$
is a Frobenius extension, this captures the notion of depth three for a Frobenius extension in \cite{KN, KS} such that if $B \| C$ is depth three, then $A \| C$
is depth two (a phenomenon of finite depth subfactors, see \cite{NV}).  We provide a similar definition
of finite depth Frobenius extension with   embedding theorem
utilizing a depth three subtower of the Jones tower. 
If $A$, $B$ and $C$ correspond to a tower of subgroups $G > H > K$ via 
the group algebra over a fixed base ring, the depth three condition
is the condition that subgroup $K$ has normal closure $K^G$ contained in $H$.
For a depth three tower of rings, there is  a pre-Galois theory for the ring $\End {}_BA_C$ and coring $(A \o_B A)^C$ involving Morita context bimodules and left coideal subrings. This is applied
in  two sections to a specialization of a   Jacobson-Bourbaki correspondence theorem
for augmented rings to  depth two extensions with depth three intermediate division rings.        
\end{abstract} 
\maketitle

\section{Introduction}

Depth two theory is a type of Galois theory for noncommutative ring
extensions, where the Galois group in field theory is replaced by a 
Hopf algebroid (with or perhaps without antipode).  
The Galois theory of depth two ring extensions has been studied
in a series of papers by the author \cite{LK2006A, LK2006B, LK2007, LK2007b}
in collaboration with Nikshych \cite{KN, KN2}, Szlach\'anyi \cite{KS},
and K\"ulshammer \cite{KK},
with a textbook treatment by Brzezi\'nski
and Wisbauer \cite{BW}. There are a number of issues that remain unexplored or unanswered in full including chirality \cite{LK2006A, LK2006B},
normality \cite{KK, LK2007}, a Galois inverse problem and  a Galois
correspondence problem \cite{Sz}.  

The Galois correspondence problem confronts the Galois theorist with a tower of three
rings $C \subseteq B \subseteq A$.  With no further assumption on the rings,
we should impose at least a relative condition on the tower to arrive at results. With this in mind we propose to generalize the notion of depth two (D2) ring extension
$A \| B$ to a notion of depth three (D3) tower $A \| B \| C$. In more
detail, the tower $A \| B \| C$ is right depth three (rD3) if $A \o_B A$
is $A$-$C$-isomorphic to a direct summand of $A \oplus \cdots \oplus A$
(finitely many times).  Many depth two
theoretic results generalize suitably, such as a decomposition of an endomorphism ring  into a crossed product with a quantum algebraic structure.  If $A \supseteq C$ is D2 and $B$ a D3 intermediate ring, we may make
use of Jacobson-Bourbaki theorems pairing certain ring extensions,
such as division rings or simple algebras, with their endomorphism rings,
in order to obtain theorems pairing D3 intermediate rings $B$ with left coideal 
subrings $\End {}_CA_B$ of the bialgebroid $\End {}_CA_C$ over the centralizer $A^C$.  

The notion of D3 tower will also serve to give a transparent and workable algebraic
definition of finite depth, originally an analytic notion in subfactor theory.
A finite Jones index subfactor may be thought of algebraically as a Frobenius extension, where the conditional expectation and Pimsner-Popa orthonormal bases are the Frobenius coordinate system.  If $B \| C$ is then a Frobenius extension with $A = \End B_C$ and $B \into A$ the left regular representation, then a depth three tower $A \| B \| C$ is a depth three
Frobenius extension (cf.\ \cite{KN} and \cite[preprint version]{KS}). 
More generally, $B \| C$ is depth $n \geq 2$ if $A_{n-2} \| A_{n-3} \| C$
is a depth three tower, where 
\begin{equation}
\label{eq: tower}
C \into B \into A \into A_2 \into \cdots \into A_n \into \cdots
\end{equation}
is the Jones tower of iterated right endomorphism rings (and where
$A_2 = \End A_B$, $\ldots$, $A = A_0$ and $B = A_{-1}$). We have the following algebraic generalization of the embedding theorem  for subfactors of Nikshych and
Vainerman \cite{NV}:  if $C \into B$ is depth $n$,
then  $C \into A_{m}$ is a depth two Frobenius extension
for some $m \geq n-2$.  

The paper is organized as follows.  
In section~2  we note that right or left D3 ring towers are characterized in terms either
of the tensor-square, H-equivalent modules, quasibases or the endomorphism ring. We prove a Theorem~\ref{th-conv} that
a depth three Frobenius extension $B \| C$
embeds in a depth two extension
$A \| C$ (where $A = \End B_C$).
In a more technical section~8 we extend
this technique to define a finite depth Frobenius extension and prove an embedding theorem for these as well: this
answers a problem raised by Nikshych and the author in \cite[Remark 5.2]{KN2}.  In section~3
we show that a tower of subgroups $G > H > K$ of finite index with the condition that the
normal closure $K^G < H$ ensures that the group algebras $ F[G] \supseteq F[H] \supseteq F[K]$
are a depth three tower w.r.t.\ any base ring $F$.  We propose that the converse
is true if $G$ is a finite group and $F = \C$.  In section~4 we study the right
coideal subring $E = \End {}_BA_C$
as well as the bimodule and co-ring $P = (A \o_B A)^C$, which provide the quasibases for a right D3 tower $A \| B \| C$.  
We show that right depth three towers may be characterized by $P$ being finite projective as
a left module over the centralizer $V = A^C$ and a pre-Galois isomorphism $A \o_B A \stackrel{\cong}{\longrightarrow} A \o_V P$.  

In section~5 we study further Galois properties of D3 towers, such as the smash product decomposition of an endomorphism ring and the invariants
as a bicommutator.  In section~6, we generalize the Jacobson-Bourbaki correspondence, which
associates $\End E_F$ to subfields
$F$ of $E$ (or skew fields), and conversely
associates $\End {}_{\mathcal{R}}E$ to
closed subrings $\mathcal{R} \subseteq \End E_F$.  We then compose this correspondence
with an anti-Galois correspondence to 
prove the main Theorem~\ref{th-GalCor}:
viz., there is a Galois correspondence
between D3 intermediate division rings
of a D2 extension of an augmented ring $A$
over a division ring $C$, on the one hand,
with Galois left coideal subrings of
the bialgebroid $\End {}_CA_C$, on
the other hand.  In Section~7, we apply
 Jacobson-Bourbaki correspondence to show that the Galois connection for separable field extensions in \cite{Sz} is a Galois correspondence between weak Hopf subalgebras and intermediate fields.

\section{Definition and first properties of depth three towers}

Let $A$, $B$ and $C$ denote rings with identity element,
and $C \to B$, $B \to A$ denote ring homomorphisms preserving the
identities.  We use ring extension notation $A \| B \| C$ for
$ C \to B \to A$ and call this a tower of rings: an important
special case if of course $C \subseteq B \subseteq A$ of subrings
$B$ in $A$ and $C$ in $B$.  Of most importance to us are the induced
bimodules such as ${}_BA_C$ and ${}_CA_B$.  We may naturally also
choose to work with algebras over commutative rings, and obtain almost
identical results.

We denote the centralizer subgroup of a ring $A$ in an $A$-$A$-bimodule $M$ by
$M^A =\{ m \in M \| \, \forall a \in A, ma = am \}$.  We also
use the notation $V_A(C) = A^C$ for the centralizer subring of $C$ in $A$.
This should not be confused with our notation $K^G$ for the normal
closure of a subgroup $K < G$.  Notation like $\End B_C$
will denote the ring of endomorphisms of the module $B_C$
under composition and  addition. We let $N^n_R$ denote the $n$-fold direct
sum of a right $R$-module $N$ with itself; let $M_R \oplus * \cong N^n_R$
denote the module $M$ is isomorphic to a direct summand of $N^n_R$. 
Finally, the symbol $\cong$ denotes  isomorphism and occasionally will
denote anti-isomorphism when we can safely
ignore opposite rings (such as ``two anti-isomorphisms compose to give an isomorphism,'' or ``two opposite rings are Morita equivalent iff the rings are Morita
equivalent'').        

\begin{definition}

A tower of rings $A \| B \| C$ is right depth three (rD3) if the tensor-square
$A \o_B A$ is isomorphic as $A$-$C$-bimodules to a direct summand of
 a finite direct
sum of $A$ with itself: in module-theoretic symbols, this becomes, for some positive integer $N$, 
\begin{equation}
\label{eq: rD3}
{}_A A\o_B A_C \oplus * \cong {}_A A^N_C
\end{equation}
\end{definition}
By switching to $C$-$A$-bimodules instead, we similarly define
a \textit{left D3 tower} of rings. The theory for these is dual to that for 
 rD3 towers; we briefly consider it at the end of this section.  As an alternative
to refering to a rD3 tower $A \| B \| C$,
we may refer to $B$ as an rD3 \textit{intermediate ring} of $A \| C$,
if $C \to A$ factors through $B \to A$ and 
$A \| B \| C$ is rD3.  

Recall that over a ring $R$, two modules $M_R$ and $N_R$ are
H-equivalent if $M_R \oplus * \cong N^n_R$ and $N_R \oplus * \cong M_R^m$
for some positive integers $n$ and $m$.  In this case, the endomorphism
rings $\End M_R$ and $\End N_R$ are Morita equivalent with
context bimodules $\Hom (M_R, N_R)$ and $\Hom (N_R, M_R)$.

\begin{lemma}
A tower $A \| B \| C$ of rings is rD3 iff the natural $A$-$C$-bimodules
$A \o_B A$ and $A$ are H-equivalent.
\end{lemma}
\begin{proof}  
We note that for any tower of rings, $A \oplus * \cong A \o_B A$ as $A$-$C$-bimodules, since the epi $\mu: A \o_B A \to A$ splits as an
$A$-$C$-bimodule arrow.  
\end{proof}

Since for any tower of rings $\End {}_AA_C$ is isomorphic to the centralizer
$V_A(C) = A^C$ (or anti-isomorphic according to convention),
we see from the lemma that the notion of rD3 has something to do with
classical depth three. Indeed,

\begin{example}
\begin{rm}
If $B \| C$ is a Frobenius extension, with Frobenius system $(E,x_i,y_i)$
satisfying for each $a \in A$,  
\begin{equation}
\sum_i E(a x_i)y_i = a = \sum_i x_i E(y_i a)
\end{equation}
then $ B \o_C B \cong \End B_C := A$ via $x \o_B y \mapsto \lambda_x \circ E \circ \lambda_y$ for left multiplication $\lambda_x$ by element $x \in B$.
Let $B \to A$ be the mapping $B \into \End B_C$ given by $b \mapsto \lambda_b$.
It is then easy to show that ${}_AB \o_C B \o_C B_C \cong {}_AA \o_B A_C$,
so that for Frobenius extensions,
 condition~(\ref{eq: rD3}) is equivalent to the condition for rD3 in preprint \cite{KS}, which in turn slightly generalizes the condition in \cite{KN}
for depth three free Frobenius extension.   
\end{rm}
\end{example}

Another litmus test for a correct notion of depth three is that depth two
extensions should be depth three
in a certain sense. Recall that
a ring extension $A \| B$ is right depth two (rD2)
if the tensor-square $A \o_B A$ is $A$-$B$-bimodule
isomorphic to a direct summand of $N$ copies of $A$ in a direct sum with itself:
\begin{equation}
\label{eq: rD2}
{}_AA \o_B A_B \oplus * \cong {}_AA^N_B
\end{equation}
 Since the notions pass from ring extension
to tower of rings, there are several cases to look at.  

\begin{prop}
\label{prop-d2d3}
Suppose $A \| B \| C$ is a tower of rings. We note:
\begin{enumerate}
\item If $B= C$ and $B \to C$ is the identity mapping,
then $A \| B \| C$ is rD3 $\Leftrightarrow$ $A \| B$ is rD2. 
\item If $A \| B$ is rD2, then $A \| B \| C$ is rD3 w.r.t.\ any ring extension
$B \| C$.
\item If $A \| C$ is rD2 and $B \| C$ is a separable extension,
then $A \| B \| C$ is rD3.    
\item If $B \|C$ is left D2, and $A = \End B_C$, then 
$A \| B \| C$ is left D3.
\item If $C$ is the trivial subring, any ring extension $A \| B$,
where ${}_BA$ is finite projective, together
with $C$ is rD3.  
\end{enumerate}
\end{prop}
\begin{proof}
The proof follows from comparing eqs.~(\ref{eq: rD3}) and~(\ref{eq: rD2}),
noting that $A \o_B A \oplus * \cong A \o_C A$ as natural $A$-$A$-bimodules
if $B \| C$ is a separable extension (thus having a separability element 
$e = e^1 \o_C e^2 \in (B \o_C B)^B$ satisfying $e^1 e^2 = 1$),
and finally from \cite{LK2006A} that $B \| C$ left D2 extension $\Rightarrow$  $A \| B$ is left D2 extension if $A = \End B_C$. The last
statement follows from tensoring ${}_BA \oplus * \cong {}_BB^n$ by
${}_AA\o_B -$.  
\end{proof}

The next theorem is a converse and algebraic
simplification of a key fact in subfactor Galois theory (the n = 3 case):  a depth three subfactor $N \subseteq M$ yields a depth
two subfactor $N \subseteq M_1$, w.r.t.\
its basic construction $M_1 \cong M \o_N M$. In preparation,
let us call a ring extension $B \| C$ rD3 if the endomorphism
ring tower $A \| B \| C$ is rD3, where $A = \End B_C$ and
$A \| B$ has underlying map $\lambda:
B \rightarrow \End B_C$, the left regular mapping
 given by $\lambda(x)(b) = xb$
for all $x, b \in B$. (This definition will
extend to identify depth $n > 3$ Frobenius
extensions as well in section~8.) 

\begin{theorem}
\label{th-conv}
Suppose $B \| C$ is a Frobenius extension
and $A = \End B_C$.  If $ B \| C$ is
rD3, then the composite extension $A \| C$ is D2. 
\end{theorem}
\begin{proof}
Begin with the well-known bimodule isomorphism for a Frobenius
extension $B \| C$,  between its endomorphism ring
and  its tensor-square, 
${}_BA_B \cong  {}_B B \otimes_C B_B$.
Tensoring by ${}_AA \o_B - \o_B A_A$,
we obtain $A \o_C A \cong A \o_B A \o_B A$
as natural $A$-$A$-bimodules.  Now
restrict the bimodule isomorphism in eq.~(\ref{eq: rD3}) on the left to $B$-modules and tensor by ${}_AA \o_B -$
to obtain ${}_AA \o_C A_C \oplus * \cong
{}_A A \o_B A_C^N$ after substitution of
the tensor-cube over $B$ by the tensor-square
over $C$.  By another application of  eq.~(\ref{eq: rD3})
we arrive at $${}_AA \o_C A_C \oplus * \cong
{}_A A^{N^2}_C$$
Thus $A \| C$ is right D2.  Since it is
a Frobenius extension as well, it is also
left depth two.  
\end{proof}

We introduce quasibases for right depth three towers.

\begin{theorem}
\label{th-qb}
A tower $A \| B \| C$ is right depth three iff there are $N$ elements each
of $\gamma_i \in \End {}_BA_C$ and of $u_i \in (A \o_B A)^C$ satisfying
(for each $x, y \in A$)
\begin{equation}
\label{eq: rd3qb}
x \o_B y = \sum_{i=1}^N x \gamma_i(y) u_i
\end{equation}
\end{theorem}
\begin{proof}
From the condition~(\ref{eq: rD3}), there are obviously $N$ maps each of 
\begin{equation}
f_i \in \Hom ({}_AA_C, {}_AA \o_B A_C), \ \ g_i \in \Hom ({}_AA \o_B A_C, {}_AA_C)
\end{equation}
such that $\sum_{i = 1}^N f_i \circ g_i = \id_{A \o_B A}$.
First, we note that for any tower of rings, not necessarily rD3, 
\begin{equation}
\Hom ({}_AA_C, {}_AA \o_B A_C) \cong (A \o_B A)^C
\end{equation}
via $f \mapsto f(1_A)$.  The inverse is given by $p \mapsto ap$
where $p = p^1 \o_B p^2 \in (A \o_B A)^C$ using a Sweedler-type notation
that suppresses a possible summation over simple tensors.

The other hom-group above also has a simplification.  We note
that for any tower,
\begin{equation}
\Hom ( {}_AA \o_B A_C, {}_AA_C) \cong \End {}_BA_C
\end{equation}
via $F \mapsto F(1_A \o_B -)$.  Given $\alpha \in \End {}_BA_C$,
we define an inverse sending $\alpha$ to the homomorphism
$x \o_B y \mapsto x \alpha(y))$.  

Let $f_i$ correspond to $u_i \in (A \o_B A)^C$ 
and $g_i$ correspond to $\gamma_i \in \End {}_BA_C$ via the mappings
just described.  We compute:
$$ x \o_B y = \sum_i f_i (g_i(x \o y)) = \sum_i f_i (x \gamma_i(y))=
\sum_i x \gamma_i(y)u_i, $$
which establishes the rD3 quasibases equation in the theorem,
given an rD3 tower.

For the converse, suppose we have $u_i \in (A \o_B A)^C$ and
$\gamma_i \in \End {}_BA_C$ satisfying the equation in the theorem.
 Then map $\pi: A^N \to A \o_B A$ by $$\pi: (a_1,\ldots,a_N) \longmapsto
\sum_i a_iu_i,$$ an $A$-$C$-bimodule epimorphism split
by the mapping $\sigma: A \o_B A \into A^N$ given by
$$\sigma(  x \o_B y ) := (x\gamma_1(y),\ldots,x\gamma_N(y)). $$
It follows from the equation above
 that $\pi \circ \sigma = \id_{A \o_B A}$.
\end{proof}

\subsection{Left D3 towers and quasibases} A tower of rings $A \| B \| C$
is left D3 (or $\ell$D3) if the tensor-square $A \o_B A$ is an $C$-$A$-bimodule direct
summand of $A^N$ for some $N$.  If $B = C$, this recovers the definition
of a left depth two extension $A \| B$.  There is a left version of
all results in this paper:  we note that  $A \| B \| C$ is a right D3
tower if and only if $A^{\rm op} \| B^{\rm op} \| C^{\rm op}$ is
a left D3 tower (cf.\ \cite{LK2006A}). 
Finally we define a tower $A \| B \| C$ to be D3 if it is both left D3 and right D3. 

The  notation in the theorem below refers to that in the example above.  
\begin{theorem}
Suppose $B \| C$ is a Frobenius extension with $A = \End B_C$.  
 Then $A \| B \| C$ is right depth three if and only if 
$A \| B \| C$ is left depth three.
\end{theorem}
\begin{proof}
It is well-known that also $A \| B$ is a Frobenius extension. 
Then $A \o_B A \cong \End A_B$ as natural $A$-$A$-bimodules.

Now note the following characterization of left D3 with proof
almost identical with that of \cite[Prop.\ 3.8]{LK2007}:
If $A \| B \| C$ is a tower where $A_B$ if finite projective, then
$A \| B \| C$ is left D3 $\Leftrightarrow$ $\End A_B \oplus * \cong A^N$
as natural $A$-$C$-bimodules.
The proof involves noting that $\End A_B \cong \Hom (A \o_B A_A,A_A)$ as natural
$A$-$C$-bimodules 
via $$f \longmapsto (a \o a' \mapsto f(a)a').$$
The finite projectivity is used for reflexivity in hom'ming this isomorphism,
thus proving the converse statement.     
 
Of course a Frobenius extension satisfies the finite projectivity condition.
Comparing the isomorphisms of $\End A_B$ and $A \o_B A$ to
direct summands of finitely many copies of $A$ just above and in
eq.~(\ref{eq: rD3}), we note that the tower is $\ell$D3 $\Leftrightarrow$
rD3.  
\end{proof}

In a fairly obvious reversal to
 opposite ring structures in the proof of Theorem~\ref{th-qb},
we see that a tower $A \| B \| C$ is left D3 iff there are
$N$ elements $\beta_j \in \End {}_CA_B$ and $N$ elements $t_j \in (A \o_B A)^C$
such that for all $x, y \in A$, we have
\begin{equation}
\label{eq: ld3qb}
x \o_B y = \sum_{j=1}^N t_j \beta_j(x)y
\end{equation}

We note explicitly that  
if $A \| B$ is a Frobenius extension with Frobenius system $(E, x_i,y_i)$, then $A \| B \| C$ is rD3 iff
the tower is $\ell$D3.  For example,
starting with the $\ell$D3 quasibases data
above, a right D3 quasibases is given by
\begin{equation}
\label{eq: left-right}
\{ E(- t^1_j)t^2_j \} \ \ \ \{ \sum_i \beta_j(x_i) \o_B y_i \} 
\end{equation}
as one may readily compute.    
 
We record the characterization of left D3, noted above in the proof, for towers satisfying a finite
projectivity condition.

\begin{theorem}
\label{th-char}
Suppose $A \| B \| C$ is a tower of rings where $A_B$ is finite projective.
Then this tower is left D3 if and only if the natural $A$-$C$-bimodules satisfy for some $N$, 
\begin{equation}
\End A_B \, \oplus \,  * \, \cong A^N
\end{equation}
\end{theorem}

Dually we establish that if $A \| B \| C$ is a tower where ${}_BA$ is finite projective,
then $A \| B \| C$ is right D3 if and only if $\End {}_BA \oplus * \cong A^N$
as natural $C$-$A$-bimodules.

%%%%%%%%%%%%%%%%%%%%%%%%%%%%%%%%%%%%%%%

\section{Depth three for towers of groups}

Fix a base ring $F$.  Groups give rise to rings via $G \mapsto
F[G]$, the functor associating the group algebra $F[G]$ to a group
$G$.  Therefore we can pull back the notion of depth $2$ or $3$
for ring extensions or towers to the category of groups when
reference is made to the base ring. 

In the paper \cite{KK}, a depth two subgroup w.r.t.\ the complex numbers
is shown to be equivalent to the notion of normal subgroup for finite
groups.  This consists of two results.  The easier result is that over any
base ring, a normal subgroup of finite index is depth two
by exhibiting left or right D2 quasibases via  coset representatives
and projection onto cosets.  This proof suggests that the converse hold
as well.  The second result is a converse for complex finite dimensional D2
group algebras where normality of the subgroup is established using character theory and Mackey's subgroup theorem.

In this section, we will similarly do the first step in showing
what group-theoretic notion corresponds to depth three tower of rings.
Let $G > H > K$ be a tower of groups, where $G$ is a finite group,
$H$ is a subgroup, and $K$ is a subgroup of $H$.  Let $A = F[G]$,
$B = F[H]$ and $C = F[K]$.  Then $A \| B \| C$ is a tower of rings,
and we may ask what group-theoretic notion on $G > H > K$ will guarantee,
with fewest possible hypotheses, that $A \| B \| C$ is rD3.

\begin{theorem}
The tower of groups algebras $A \| B \| C$ is D3 if
the corresponding tower of groups $G > H > K$ satisfies
\begin{equation}
K^G < H
\end{equation}
where $K^G$ denotes the normal closure of $K$ in $G$.  
\end{theorem}
\begin{proof}
Let $\{ g_1,\ldots,g_N \}$ be double coset representatives such
that $G = \coprod_{i=1}^N Hg_iK$.  Define $\gamma_i(g) = 0$
if $g \not \in Hg_i K$ and $\gamma_i (g) = g$ if $g \in Hg_iK$.  
Of course, $\gamma_i \in \End {}_BA_C$ for $i = 1,\ldots, N$.  

Since $K^G \subseteq H$, we have $gK \subseteq Hg$ for each $g \in G$.
  Hence for each $k \in K$, $g_jk = h g_j$ for some $h \in H$.  It follows
that $$g_j^{-1} \o_B g_jk = g_j^{-1} h \o_B g_j = kg_j^{-1} \o_B g_j.  $$

Given $g \in G$, we have $g = hg_j k$ for some $j = 1,\ldots,N$,
$h \in H$, and $k \in K$.
Then we compute:
$$ 1 \o_B g = 1 \o_B hg_jk = hg_j g_j^{-1} \o_B g_jk = hg_jk g_j^{-1} \o_B g_j $$
so $1 \o_B g =  \sum_i \gamma_i(g) g_i^{-1} \o_B g_i $
where $g_i^{-1} \o_B g_i \in (A \o_B A)^C$.  By theorem then,
$A \| B \| C$ is an rD3 tower.

The proof that the tower of group algebras is left D3 is entirely symmetical
via the inverse mapping.    
\end{proof}

The theorem is also valid for infinite groups where the index $[G : H]$
is finite, since $HgK = Hg$ for each $g \in G$.  

Notice how the equivalent notions of depth two and normality
for finite groups over $\C$ yields the Proposition~\ref{prop-d2d3}
for groups.  Suppose we have a tower of groups $G > H > K$ where
$K^G \subseteq H$.  If $K = H$, then $H$ is normal (D2) in $G$.
If $K = \{ e \}$, then it is rD3 together with any subgroup $H < G$.
If $H \ract G$ is a normal subgroup, then necessarily $K^G \subseteq H$.  
If $K \ract G$, then $K^G = K < H$ and the tower is D3.

Question:  Can the character-theoretic proof in \cite{KK} be adapted
to prove that a D3 tower $\C[G] \supseteq \C[H] \supseteq \C[K]$
where $G$ is a finite group 
satisfies $K^G < H$?

%%%%%%%%%%%%%%%%%%%%%%%%%%%%%%%%%%%%%%%%%%%%%%%%%%%%%%%%%%%%%%%%%%%%%%%%%%

\section{Algebraic structure on $\End {}_BA_C$ and $(A \o_B A)^C$}

In this section, we study the calculus of some structures definable for an
rD3 tower $A \| B \| C$, which reduce to the dual bialgebroids over
the centralizer of a ring extension in case $B = C$ and their actions/coactions.
Throughout the section, $A \| B \| C$ will denote a right depth three
tower of rings, 
$$ P := (A \o_B A)^C, \ \ Q := ( A \o_C A)^B,  $$
which are bimodules with respect to the two rings familiar in depth
two theory,
$$ T := (A \o_B A)^B, \ \ U := (A \o_C A)^C$$
Note that $P$ and $Q$ are isomorphic to two $A$-$A$-bimodule Hom-groups: 
\begin{equation}
\label{eqs: tys}
 P \cong \Hom (A \o_C A, A \o_B A), \ \ Q \cong \Hom (A \o_B A, A \o_C A).
\end{equation}
Recall that $T$ and $U$ have multiplications given by
$$ tt' = {t'}^1 t^1 \o_B t^2 {t'}^2, \ \ uu' = {u'}^1u^1 \o_C u^2 {u'}^2, $$
where $1_T = 1_A \o 1_A$ and a similar expression for $1_U$.  
Namely, the bimodule ${}_TP_U$ is given by
\begin{equation}
 {}_TP_U:\ t \cdot p \cdot u = u^1 p^1 t^1 \o_B t^2 p^2 u^2 
\end{equation}
The bimodule ${}_UQ_T$ is given by
\begin{equation}
{}_UQ_T: \ u \cdot q \cdot t = t^1 q^1 u^1 \o_C u^2 q^2 t^2
\end{equation}

We have the following result, also mentioned in passing
in \cite{LK2006B} with several additional hypotheses.

\begin{prop}
\label{prop-ME}
The bimodules $P$ and $Q$ over the rings $T$ and $U$ form
a Morita context with associative multiplications
\begin{equation}
P \o_U Q \rightarrow T,\ \ p \o q \mapsto pq = q^1 p^1 \o_B p^2 q^2 
\end{equation}
\begin{equation}
Q \o_T P \rightarrow U, \ \ q \o p \mapsto qp = p^1 q^1 \o_C q^2 p^2
\end{equation}
If $B \| C$ is an
H-separable extension, then $T$ and $U$ are Morita equivalent rings via this context.
\end{prop}
\begin{proof}
The equations $p(qp') = (pq)p'$ and $q(pq') = (qp)q'$ for
$p,p' \in P$ and $q,q' \in Q$ follow from the four equations directly above.  

Note that
$$ T \cong \End {}_AA \o_B A_A, \ \ U \cong \End {}_AA \o_C A_A $$
as rings.  We now claim that the hypotheses on $A \| B$, $A \| C$ and $B \| C$
imply that the $A$-$A$-bimodules $A \o_B A$ and
$A \o_C A$ are H-equivalent.  Then the endomorphism rings above are Morita
equivalent via context bimodules given by eqs.~(\ref{eqs: tys}), which
proves the proposition.

Since $B \| C$ is H-separable, it is in particular separable, and
the canonical $A$-$A$-epi $A \o_C A \rightarrow A \o_B A$ splits 
via an application of a separability element.  Thus, $A \o_B A \oplus * \cong
A \o_C A$.  
The defining condition for H-separability is $B \o_C B \oplus * \cong B^N$ as $B$-$B$-bimodules for some positive integer $N$. Therefore, $A \o_C A \oplus * \cong A \o_B A^N$ as $A$-$A$-bimodules
by an application of the functor $A \o_B \, - \, \o_B A$.  Hence,
$A \o_B A$ and $A \o_C A$ are H-equivalent $A^e$-modules (i.e., $A$-$A$-bimodules). 
\end{proof}
   
We denote the centralizer subrings $A^B$ and $A^C$ of $A$ by
\begin{equation}
 R := V_A(B) \subseteq V_A(C) := V
\end{equation}

From $R \cong \Hom (A \o_B A, A)$ and $V \cong \Hom (A \o_C A, A)$
and composition with eq.~(\ref{eqs: tys}), we obtain 
the generalized anchor mappings (cf.\ \cite{LK2006B}),
\begin{equation}
R \o_T P \longrightarrow V, \ \ r \o p \longmapsto p^1rp^2
\end{equation}
\begin{equation}
V \o_U Q \longrightarrow R, \ \ v \o q \longmapsto q^1 v q^2
\end{equation}
 
\begin{prop}
\label{prop-b}
The two generalized anchor mappings are bijective if
$B \| C$ is H-separable.
\end{prop}
\begin{proof}
Denote $r \cdot p := p^1 r p^2$ and $v \cdot q := q^1 v q^2$.
From the previous propostion, there are elements $p_i \in P$
and $q_i \in Q$ such that $\sum_i p_i q_i = 1_T$;
in addition, ${p'}_j \in P$ and ${q'}_j \in Q$ such that
$1_U = \sum_j {q'}_j {p'}_j $.  
Let $ v \in V$, then
$$ v = v \cdot 1_U = \sum_j v \cdot ({q'}_j {p'}_j) = \sum_j (v \cdot {q'}_j)\cdot {p'}_j $$
and a similar computation starting with $r = r \cdot 1_T$ shows that the two generalized
anchor mappings are surjective.  

In general, we have the corestriction of the inclusion $T \subseteq A \o_BA$,
 \begin{equation}
{}_TT \into {}_TP
\end{equation}
which is split as a left $T$-module monic by $p \mapsto e^1pe^2$ in case there is a separability element
$e = e^1 \o_C e^2 \in B \o_C B$.  Similarly, 
\begin{equation}
{}_UQ \into {}_UU
\end{equation}
is a split monic in case $B \| C$ is separable. If
$\sum_i v_i \o_U q_i \in \ker (V \o_U Q \rightarrow R)$
then $\sum_i v_i \o_U q_i \mapsto \sum_i v_i \cdot q_i \o_U 1_U = 0$
via an injective mapping, whence $\ker (V \o_U Q \rightarrow R) = \{ 0 \}$.  
 
Of course, if $B \| C$ is H-separable, we note from  Proposition~\ref{prop-ME} and Morita theory 
that $P$ and $Q$ are projective generators on both sides,
(and faithfully flat). If $K := \ker (R \o_T P \rightarrow V)$,
then $K \o_U Q = 0$,
since $\sum_j r_j \cdot p_j = 0$ implies
$$ \sum_j r_j \o_T p_j \o_U q \longmapsto \sum_j r_j \o_T p_j q \o_U 1_U = 0$$
via an injective mapping.  
It follows from faithful flatness of ${}_UQ$ that $K = \{ 0\}$. 
\end{proof}

Note that $P$ is a $V$-$V$-bimodule (via the commuting homomorphism
and antihomomorphism $V \rightarrow U \leftarrow V$):
\begin{equation}
{}_VP_V:\ \ v \cdot p \cdot v' = vp^1 \o_B p^2 v'
\end{equation}
Note too that $E = \End {}_BA_C$ is an $R$-$V$-bimodule via
\begin{equation}
{}_RE_V: \ \ r \cdot \alpha \cdot v = r\alpha(-)v
\end{equation}
Note the subring and over-ring
\begin{equation}
\End {}_BA_B \subseteq E \subseteq \End {}_CA_C
\end{equation}
which are the total algebras of the left $R$- and $V$- bialgebroids in depth two theory \cite{KS, LK2006A, LK2006B}.  

\begin{lemma}
\label{lemma-agema}
The modules ${}_VP$ and $E_V$ are finitely generated projective. In case $A \| C$ is
left D2, the subring
$E$ is a right coideal subring of the left $V$-bialgebroid
$\End {}_CA_C$.   
\end{lemma}
\begin{proof}
This follows from eq.~(\ref{eq: rd3qb}), since $p \in P \subseteq A \o_B A$, so 
$$ p = \sum_i p^1 \gamma_i(p^2)u_i $$
where $u_i \in P$ and $p \mapsto p^1 \gamma_i(p^2)$ is in $\Hom ({}_VP, {}_VV)$,
thus dual bases for a finite projective module.  The second claim follows similarly
from
$$ \alpha = \sum_i \gamma_i(-)u_i^1 \alpha(u^2) $$
where $\gamma_i \in E$ and $\alpha \mapsto u^1 \alpha(u^2)$
are mappings in $\Hom (E_V, V_V)$.

Now suppose $\beta_j \in S := \End {}_CA_C$
and $t_j \in (A \o_C A)^C$ are left
D2 quasibases of $A \| C$.  Recall
that the coproduct $\cop: S \rightarrow
S \o_V S$ given by ($\beta \in S$)
\begin{equation}
\cop(\beta) = \sum_j \beta(-t^1_j)t^2_j \o_V \beta_j
\end{equation}
makes $S$ a left $V$-bialgebroid \cite{KS}.
Of course this restricts and corestricts to $\alpha \in E$
as follows:  $\cop(\alpha) \in E \o_V S$.
Hence, $E$ is a right coideal subring of $S$.     
\end{proof}

In fact, if $A \| B$ is also D2,
and $\mathcal{S} = \End {}_BA_B$, 
then $E$ is similarly shown to be 
an $\mathcal{S}$-$S$-bicomodule ring
For we recall the coaction $E \rightarrow
\mathcal{S} \o_R E$ given by 
\begin{equation}
\alpha\-1 \o_R \alpha\0 = \sum_i \tilde{\gamma}_i \o \tilde{u}_i^1 \alpha(\tilde{u}_i^2 -) 
\end{equation}
where $\tilde{\gamma}_i \in \mathcal{S}$
and $\tilde{u}_i \in (A \o_B A)^B$
are right D2 quasibases of $A \| B$
(restriction of \cite[eq.~(19)]{LK2006A}).  
 
Twice above we made use of a $V$-bilinear pairing $P \o E \to V$ given
by
\begin{equation}
\bra p, \alpha \ket := p^1 \alpha(p^2), \ \ (p \in P= (A \o_B A)^C,\, \alpha \in E= \End {}_BA_C)
\end{equation}

\begin{lemma}
The pairing above is nondegenerate.
It induces $E_V \cong \Hom ({}_VP, {}_VV)$ via $\alpha \mapsto \bra - , \alpha \ket$.
\end{lemma}
\begin{proof}
The mapping has the inverse $F \mapsto \sum_i \gamma_i(-)F(u_i)$ where
$\gamma_i \in E, u_i \in P$ are rD3 quasibases for $A \| B \| C$.  
Indeed, $\sum_i \bra p, \gamma_i \ket F(u_i) =$ $ F(\sum_i p^1 \gamma_i(p^2)u_i) = F(p)$
for each $p \in P$ since $F$ is left $V$-linear,
and for each $\alpha \in E$, we note that $\sum_i \gamma_i(-)\bra u_i, \alpha\ket = \alpha$.
\end{proof}

\begin{prop}
There is a $V$-coring structure on $P$ left dual to the ring structure on $E$.
\end{prop}
\begin{proof}
We note that 
\begin{equation}
P \o_V P \cong (A \o_B A \o_B A)^C 
\end{equation}
via $p \o p' \mapsto p^1 \o p^2 {p'}^1 \o {p'}^2$
with inverse $$p = p^1 \o p^2 \o p^3 \mapsto \sum_i (p^1 \o_B p^2 \gamma_i(p^3)) \o_V u_i.$$

Via this identification, define a $V$-linear coproduct $\cop: P \rightarrow P \o_V P$
by \begin{equation}
\cop(p) = p^1 \o_B 1_A \o_B p^2.
\end{equation}
Alternatively, using Sweedler notation and rD3 quasibases,
\begin{equation}
\label{eq: cop}
p\1 \o_V p\2 = \sum_i (p^1 \o_B \gamma_i(p^2)) \o_V u_i 
\end{equation}
   Define a $V$-linear counit $\eps: P \to V$
by $\eps(p) = p^1p^2$.  The counital equations follow readily \cite{BW}.  

Recall from Sweedler \cite{Sw} that the $V$-coring $(P,V,\cop, \eps)$ has
left dual ring ${}^*P := \Hom ({}_VP, {}_VV)$ given by Sweedler notation 
by 
\begin{equation}
(f * g)(p) = f(p\1 g(p\2))
\end{equation}
with $1 = \eps$. Let $\alpha, \beta \in E$.
 If $f = \bra -, \alpha \ket$
and $g = \bra - , \beta \ket$ , we compute $f * g = \bra - , \alpha \circ \beta \ket$
below, which verifies the claim:  
$$ f(p\1 g(p\2)) = \sum_i \bra p^1 \o_B \gamma_i(p^2) \bra u_i, \beta \ket,\, \alpha \ket = 
\bra p^1 \o_B \beta(p^2), \alpha \ket = \bra p, \alpha \circ \beta \ket . \qed $$
\renewcommand{\qed}{}\end{proof}

In addition, we note that $P$ is $V$-coring with grouplike element
\begin{equation}
g_P := 1_A \o_B 1_A
\end{equation}
since $\cop(g_P) = 1 \o 1 \o 1 = g_P \o_V g_P$ and $\eps(g_P) = 1$.
  
There is a pre-Galois structure on $A$ given by the right $P$-comodule
structure $\delta: A \rightarrow A \o_V P$, $\delta(a) = a\0 \o_V a\1$
defined by 
\begin{equation}
\delta(a):= \sum_i \gamma_i(a) \o_V u_i.
\end{equation}
The pre-Galois isomorphism $\beta: A \o_BA \stackrel{\cong}{\longrightarrow}
 A \o_V P$ given by \begin{equation}
\beta(a \o a') = a {a'}\0 \o_V {a'}\1
\end{equation}
 is utilized below 
in another characterization of right depth three towers.

\begin{theorem}
A tower of rings $A \| B \| C$ is right depth three if and only if
${}_VP$ is finite projective and $A \o_V P \cong A\o_B A$ as natural
$A$-$C$-bimodules. 
\end{theorem}
\begin{proof}
($\Leftarrow$) If ${}_VP \oplus * \cong {}_VV^N$ and $A \o_V P \cong A \o_B A$,
then tensoring by $A \o_V -$, we obtain $A \o_B A \oplus * \cong A^N$ as natural
$A$-$C$-bimodules, the rD3 defining condition on a tower.

($\Rightarrow$) By lemma  ${}_VP$ is 
f.g.\ projective. Map $A \o_V P \rightarrow A \o_B A$ by $a \o p \mapsto
ap^1 \o_B p^2$, clearly an $A$-$C$-bimodule homomorphism.  
The inverse is the ``pre-Galois'' isomorphism, 
\begin{equation}
\beta: \, A \o_B A \rightarrow A \o_V P, \ \ \beta(a \o_B a') = \sum_i a \gamma_i(a') \o_V  u_i
\end{equation}
 since $\sum_i a p^1 \gamma_i(p^2) \o_V u_i = a \o_V p$ and
$\sum_i a \gamma_i(a') u_i = a \o a'$ for $a,a' \in A, p \in P$.  
\end{proof}

If $B \| C$ is H-separable, there is more to say about the structure of
the Morita equivalent total rings for the bialgebroids $T$ and $U$ and bijective anchor maps in propositions~\ref{prop-ME} and~\ref{prop-b}.  This stems from the fact
that for a $B$-bimodule $M$, we have an Azumaya-type
 condition for the centralizers,  $M^C \cong M^B \otimes_{Z(B)} B^C$ via $m \o c \mapsto mc$ in one direction.  This may now be applied
to each of the cases $M = A$, $A \o_B A$, and $A \o_C A$ to obtain formulas
relating $V$ and $R$, $T$ and $P$, as well as $Q$ and $U$.  We will
study the relationship of these remarks to monoidal functors and Takeuchi's $\sqrt{\rm
Morita}$ base change outlined in the paper \cite{Sch} in another paper.  
 
%%%%%%%%%%%%%%%%%%%%%%%%%%%%%%%%%%%%%%

\section{Further Galois properties of
depth three} 

We will show here that the smaller of the endomorphism rings of a depth three tower decomposes tensorially
over the overalgebra and the mixed bimodule
endomorphism ring studied above.  In case
the composite ring extension is depth
two, this is a smash product decomposition
in terms of a coideal subring of
a bialgebroid.  Finally, we express
the invariants of this coideal subring
acting on the overalgebra in terms of a
bicommutator.  

\begin{theorem}
\label{th-end}
If $A \| B \| C$ is left D3, then
\begin{equation}
\End A_B \cong A \o_V \End {}_CA_B
\end{equation}
via the homomorphism $A \o_V \End {}_CA_B
\rightarrow \End A_B$ given by $a \o_V \alpha \mapsto \lambda_a \circ \alpha$.
\end{theorem}
\begin{proof}
Given a left D3 quasibases $\beta_j \in \End {}_CA_B$ and $t_j \in (A \o_B A)^C$,
note that the mapping $\End A_B \rightarrow
A \o_V \End {}_CA_B$ given by
\begin{equation}
f \longmapsto \sum_j f(t^1_j)t^2_j \o_V \beta_j
\end{equation}
is an inverse to the homomorphism above.   
\end{proof}

\begin{cor}
\label{cor-endosmash}
If $A \| C$ is additionally D2, then
$\End {}_CA_B$ a left coideal subring
of $\End {}_CA_C$ and there is a
ring isomorphism with a smash product ring,
\begin{equation}
\End A_B \cong A \rtimes \End {}_CA_B
\end{equation}
\end{cor}
\begin{proof}
Recall from depth two theory \cite{KS} that the
$V$-bialgebroid $\End {}_CA_C$ acts
on the module algebra $A$ by simple
evaluation, $\beta \lact a = \beta(a)$.
That the action is measuring is not
hard to see from the formula for the coproduct on $\End {}_CA_C$ given
by 
\begin{equation}
\cop(\beta) = \beta\1 \o_V \beta\2 := 
\sum_k \tilde{\gamma}_k \o_V \tilde{u}^1_k
\beta(\tilde{u}^2_k-) 
\end{equation}
where $\tilde{\gamma}_k \in \End {}_CA_C$
and $\tilde{u}_k \in (A \o_C A)^C$
are right D2 quasibases for the
composite ring extension $A \| C$.
Note then that for $\alpha \in \End {}_CA_B 
\subseteq \End {}_CA_C$,
the equation yields $\alpha\1 \o_V \alpha\2
\in \End {}_CA_C \o_V \End {}_CA_B$.
Hence, $\End {}_CA_B$ is a left coideal
subring.  The details and verifications of the definition
of such an object, over a smaller base ring
than that of the bialgebroid, are rather straightforward and left to the reader.

As a consequence of the smash product
formula $\End A_C \cong A \rtimes \End {}_CA_C$ over the centralizer $V$, we
restrict to $\End A_B \subseteq \End A_C$,
apply the theorem above, to obtain
the equation for $\alpha, \beta \in \End {}_CA_B$,
\begin{equation}
(a \# \alpha)(b \# \beta) = a(\alpha\1 \lact b) \# \alpha\2 \circ \beta \in
A \o_V \End {}_CA_B
\end{equation}
where $a, b \in A$, and $\rtimes$,  $\#$ are used interchangeably.  
\end{proof}

In case $A \| C$ continues to be a D2
extension, the theorem below will characterize the subring $A^S$ of
invariants of $S = \End {}_CA_C$
as well as $A^{\mathcal{J}}$ where
$\mathcal{J} := \End {}_CA_B$,
the coideal subring of $S$, in terms
of $A$ as the natural module over $E :=
\End A_B$. The endomorphism ring $\End {}_EA$ is familiar  from the Jacobson-Bourbaki theorem in Galois theory \cite{J, P}. 

\begin{theorem}
\label{th-bic}
Let $A \| B \| C$ be left D3 and
$$A^{\mathcal{J}} = \{ x \in A | \forall
\alpha \in \mathcal{J}, \alpha(x) = \alpha(1) x \}.$$
Then $A^{\mathcal{J}} \cong \End {}_EA$
via the anti-isomorphism $x \mapsto \rho_x$.
\end{theorem}
\begin{proof}
We first note that $A^{\mathcal{J}} =
\{ x \in A | \forall f \in E, y \in A,
f(yx) = f(y)x \}$.
The inclusion $\supseteq$ easily follows
from letting $y = 1_A$ and $\alpha \in \mathcal{J} \subseteq E$.  The reverse
inclusion follows from Theorem~\ref{th-end}. Since $E \cong
A \o_V \mathcal{J}$, note that $f \circ \lambda_y \in E$ decomposes as $\sum_j f(yt^1_j)t^2_J \o \beta_j \in A \o_V \mathcal{J}$ for an arbitrary
$y \in A$.  Given $x \in A$ such that  $\alpha(x) = \alpha(1)x$
for each $\alpha \in \mathcal{J}$,
then $$f(yx) = \sum_j f(yt^1_j)t^2_j \beta_j(x) = \sum_j f(yt^1_j)t^2_j \beta_j(1)x =
f(y)x.$$

It follows from these considerations
that $\rho_x \in \End {}_EA$ for
$x \in A^{\mathcal{J}}$, since $\rho_x(f(a))
= f(\rho_x(a))$ for each $f \in E, a \in A$.

Now an inverse mapping $\End {}_EA \rightarrow A^{\mathcal{J}}$ is given
by $G \mapsto G(1)$.  Of course $\rho_x(1) = x$.  
Note that $G(1) \in A^{\mathcal{J}}$,  since for
 $\alpha \in \mathcal{J}$,
 we have $\alpha(G(1)) = G(\alpha(1))
=\lambda_{\alpha(1)}G(1)$,
since $\lambda_a \in E$ for all $a \in A$.
Finally, we note that  $G(a) = G \circ \lambda_a (1) = a G(1)$,
 whence $G
= \rho_{G(1)}$ for each $G \in \End {}_EA$.
\end{proof}
  
The following clarifies and extends part of \cite[4.1]{KS}.  Let $S$ denote the bialgebroid
$\End {}_BA_B$ below and $E$ as before
is $\End A_B$.  

\begin{cor}
If $A \| B$ is left D2, then
$A^S \cong \End {}_EA$.  
Thus if $A_B$ is balanced, $A^S = B$.
\end{cor}
\begin{proof}
Follows by Prop.~\ref{prop-d2d3} and from the theorem by letting $B = C$.  We note additionally from its proof
that 
\begin{equation}
A^S = \{ x \in A | \forall \alpha \in S, \alpha(x) = x \alpha(1) \} 
\end{equation}
since $\rho_{\alpha(1)} \in E$ in this
case.  

If $A_B$ is balanced, $\End {}_EA = \rho(B)$
by definition. This recovers the result
in \cite[Section 4]{KS}.
\end{proof} 
 
In other words, this corollary states that
the invariant subring of $A$ under the action of the bialgebroid $S$ is (anti-isomorphic to) the bicommutator
of the natural module $A$. Sugano
studies the derived ring extension
$A^* \| B^*$ of bicommutants of a ring
extension $A \| B$, where
$M_A$ is a faithful module, $E :=
\End M_A$, $\mathcal{E} := \End M_B$, 
$A^* = \End {}_EM$, $B^* = \End {}_{\mathcal{E}}M$ and there are
natural monomorphisms $A \to A^*$
and $B \to B^*$ commuting with the
mappings $B \to A$ and $B^* \to A^*$
 \cite{Su}: in these terms, $A^S \subseteq A$ is then the bicommutator of $A_A$
over the depth two extension $A \| B$.

%%%%%%%%%%%%%%%%%%%%%%%%%%%%%%%%%%%%%%%
\section{A Jacobson-Bourbaki Correspondence
for Augmented Rings}

The Jacobson-Bourbaki correspondence is
usually given between subfields $F$ of finite codimension in a field $E$
on the one hand, and their linear endomorphism rings $\End E_F$
on the other hand.
A subring of $\End E_F$ which is itself
an endomorphism ring of this form is characterized by containing $\lambda(E)$
and being finite dimensional over this.
The inverse correspondence  associates
to such a subring $R \subseteq \End E_F$,
the subfield $\End {}_RE$, since
${}_RE$ is simple as a module.
(The centralizer or commutant of $R$ in $\End E_{\Z}$
in other words.)  The correspondences are inverse to one another
by the Jacobson-Chevalley density theorem,
and may be extended to division rings
\cite[Section 8.2]{J}.
 
Usual Galois theory follows from this
correspondence, for if $E^G = F$ where $G$ is a finite group of automorphisms of $E$, then
$\End E_F \cong E \# G$ and subrings
of the form $\End E_K$ correspond
to the subrings $E \# H$ where $H$ is
a subgroup of $G$ such that $E^H = K$
for an intermediate field $K$ of $F \subseteq E$.  In this section, we will use a similar idea to pass from
the Jacobson-Bourbaki correspondence
to the correspondence $A \| B \mapsto \End {}_BA_B$ and inverse $S \mapsto A^S$
for certain Hopf subalgebroids $S$ of 
$\End {}_BA_B$ for certain depth two extensions $A \| B$.  First, we will
give an appropriate generalization of
the Jacobson-Bourbaki correspondence to
noncommutative algebra, with a proof
similar to Winter \cite[Section 2]{W}.

For the purposes below, we say an \textit{augmented
ring} $(A,D)$ is a ring $A$ with a ring homomorphism
$A \to D$ where $D$ is a division ring.   Examples are division rings, local rings, Hopf algebras and augmented algebras.  
A subring $\mathcal{R}$ of $\End A := \End A_{\Z}$ containing $\lambda(A)$,
left finitely generated over this,
where ${}_{\mathcal{R}}A$ is simple, 
is said to be a \textit{Galois subring}.   
 
\begin{theorem}[Jacobson-Bourbaki correspondence for noncommutative augmented rings]
\label{th-JBT}
Let $(A, D)$ be an augmented ring.
There is a one-to-one correspondence
between the set of division rings $B$ within $A$, 
where $B$ is a subring of $A$ and $A_B$
is a finite dimensional right vector space,
and the set of Galois subrings of $\End A$.  The correspondence is given 
by $B \mapsto \End A_B$ with inverse
correspondence $\mathcal{R} \mapsto \End {}_{\mathcal{R}}A$.  
\end{theorem}
\begin{proof}
We first show that if $B$ is a division
ring and subring of $A$ of finite right
codimension, then $E = \End A_B$ is a Galois subring and $\End {}_EA \cong B$.  We will need a theory
of left or (dually) right vector spaces over a division ring as for example to be found in \cite[chap. 4]{H}.
Suppose $[A: B]_r = d$.  

Since $\End A_B$ is isomorphic to square
matrices of order $d$ over the division ring $B$, it follows that $\End A_B$ is
finitely generated over the algebra $\lambda(A)$ of left multiplications
of $A$. Also ${}_EA$ is simple, since
$E = \End A_B$ acts transitively on
$A$.  Hence $\End {}_EA$ is a division
ring.  Since \begin{equation}
\label{eq: split}
A_B = B_B \oplus W_B
\end{equation}
for some complementary subspace $W$ over $B$,
it follows from Morita's lemma
(``generator modules are balanced'') that in fact $B \cong \End {}_EA$.

    Conversely, let $\mathcal{R}$ be a Galois subring.  Let $F^{\rm op} = \End {}_{\mathcal{R}}A$ be the
division ring (by Schur's lemma)
contained in $A^{\rm op}$ (since $A \subseteq \mathcal{R}$ and  $\End {}_AA \cong A^{\rm op}$). 
To finish the proof we need to show that $[A: F]_r < \infty$ and $\mathcal{R} = \End A_F$.   

Since $\mathcal{R}$ is finitely generated over $A$,
we have $s_1, \ldots s_n \in \mathcal{R}$ such 
that $$\mathcal{R} = As_1 + \cdots + As_n. $$
Let $e_1, \ldots, e_m \in A$ be linearly
independent in the right vector space
$A$ over $F$.  Since ${}_{\mathcal{R}}A$ is simple,
the Jacobson-Chevalley density theorem ensures
the existence of elements $r_1, \ldots, r_m
\in \mathcal{R}$ such that for all $i$ and $k$, 
$$ r_i(e_k) = \delta_{ik} 1_A.$$

By the lemma below and the hypothesis
that $A$ is an augmented ring, $m \leq n$.
With a maximal linear independent set
of vectors $e_i$ in $A$, we may assume
$e_1,\ldots,e_m$ a basis for $A_F$. 
By definition of $F$, we have $\mathcal{R} \subseteq \End A_F$.   
Let $E_{ij} := e_i r_j$ for $1 \leq i,j \leq m$ in $\mathcal{R} $. 
Since $E_{ij} (e_k) = \delta_{jk}e_i$,
these are matrix units which span
$\End A_F$. Hence $\End A_F = \mathcal{R}$.     
\end{proof}

\begin{lemma}
Let $s_1,\ldots,s_n \in \End A_{\Z}$ where
$(A,D)$ is an augmented ring.
Suppose that $$ r_1,\ldots, r_m \in As_1 + \cdots + As_n$$
and there are elements $e_1,\ldots, e_m \in A$ such that $r_i( e_k) = \delta_{ik}1_A$
for $1 \leq i,k \leq m$.  Then $m \leq n$.
\end{lemma}
\begin{proof}
By the hypothesis, there are elements
$a_{ij} \in A$ such that $r_i = \sum_{j=1}^n a_{ij}s_j$ for each $i = 1,\ldots,m$.  Then for $1 \leq i,k \leq m$, 
$$ \sum_{j=1}^n a_{ij} s_j(e_k) = r_ie_k = \delta_{ik} 1_A. $$
Applying the ring homomorphism $A \rightarrow D$ into the division ring $D$, where $a_{ij} \mapsto
d_{ij}$, $s_j(e_k) \mapsto z_{jk}$, we obtain
the matrix product equation, 
$$ \left( \begin{array}{ccc}
d_{11} & \cdots & d_{1n} \\
 \vdots &  \vdots & \vdots \\
d_{m1} & \cdots & d_{mn} 
\end{array} \right)
\left( \begin{array}{ccc}
z_{11} & \cdots & z_{1m} \\
 \vdots &  \vdots & \vdots \\
z_{n1} & \cdots & z_{nm} 
\end{array} \right) = 
\left( \begin{array}{ccc}
1_D & \cdots & 0 \\
 \vdots &  \ddots & \vdots \\
0 & \cdots & 1_D 
\end{array} \right)
$$
This shows in several ways that $m \leq n$;
for example, by the rank + nullity theorem
for right vector spaces \cite[Ch.\ 4, corollary 2.4]{H}.
\end{proof}

Let $A \supseteq C$ be a D2 ring extension, so that $S := \End {}_CA_C$ is canonically a left bialgebroid over the centralizer $A^C$.
Any D2 subextension $A \supseteq B$ has sub-$R$-bialgebroid $
\mathcal{S} := \End {}_BA_B$
where $R = A^B \subseteq A^C$.  If all extensions are balanced,
as in the situation we consider above, we recover the intermediate
D2 subring $B$ by $\mathcal{S} \leadsto A^{\mathcal{S}} = B$.  Whence
$B \leadsto \mathcal{S}$ is a surjective correspondence
and Galois connection between the set of intermediate D2 subrings
of $A \supseteq C$ and the set of sub-$R$-bialgebroids of $S$
where $R$ is a subring of $A^C$.  We  widen
our perspective to include D3 intermediate subrings $B$,
i.e. D3 towers $A \supseteq B \supseteq C$, and left coideal
subrings of $S$ in order to pass from surjective Galois connection
to Galois correspondence.

The Galois correspondence given by $B \leadsto \End {}_CA_B$
and $\mathcal{J} \leadsto A^{\mathcal{J}}$  factors through
the Jacobson-Bourbaki correspondence sketched in the theorem above.     
We  apply Theorems~\ref{th-bic},~\ref{cor-endosmash}
and~\ref{th-char} below to do this.  We  need
a notion of \textit{Galois left coideal subring} $\mathcal{J}$
of a left $V$-bialgebroid $S$.  For this we require of the left
coideal subring $\mathcal{J} \subseteq \End {}_CA_C$
that 
\begin{enumerate}
\item the module ${}_V\mathcal{J}$ is finitely generated projective where $V = A^C$;
\item $A$ has no proper $\mathcal{J}$-stable left ideals. 
\end{enumerate}
 
\begin{theorem}
\label{th-GalCor}
Let $A \supseteq C$ be a D2 extension of
an augmented ring $A$ over a division ring $C$, with  centralizer
$A^C$ denoted by $V$ and  left $V$-bialgebroid $\End {}_CA_C$
by $S$. Suppose $A_V$ is faithfully flat. Then the left D3 intermediate division rings of $A \supseteq C$ are
in Galois correspondence with the Galois left coideal subrings
of $S$.  
\end{theorem}
\begin{proof}
Since $C \subseteq A$ is D2 and left or right split (as in eq.~\ref{eq: split}), we may apply  a projection ${}_CA \rightarrow {}_CC$ to
the left D2 quasibases eq.\  
to see that $A_C$ is a finite
dimensional right vector space. For the same reasons,
each extension $A \supseteq B$ (for an intermediate
division ring $B$) is balanced by Morita's lemma.
If $B$ is additionally a left D3 intermediate ring, with
 $\mathcal{J} = \End {}_CA_B$ a left coideal subring of the
bialgebroid $S$ by Corollary~\ref{cor-endosmash} 
we have by
Theorem~\ref{th-bic} that the invariant subring
$A^{\mathcal{J}} = B$. We just note that ${}_V\mathcal{J}$ is
f.g.\ projective by the  dual of Lemma~\ref{lemma-agema},
and that a proper $\mathcal{J}$-stable left ideal of $A$
would be a proper $\End A_B$-stable left ideal in contradiction
of the transitivity argument in Theorem~\ref{th-JBT}. Thus
$B \mapsto \End {}_CA_B$ is a surjective order-reversing correspondence
between the set of left D3 intermediate division rings $A \supseteq B \supseteq C$ into the set of Galois
left coideal subrings of the $V$-bialgebroid $S$.

Suppose we are given a Galois left coideal subring $\mathcal{I}$
of $S = \End {}_CA_C$.  Then the smash product ring $A \rtimes \mathcal{I}$
has image we denote by $\mathcal{R}$ in $\End A_C$ via $a \o_V \alpha \mapsto \lambda_a \circ \alpha$ that is clearly a Galois subring, since $\lambda(A) \subseteq \mathcal{R}$
and is a finitely generated extension; also the module ${}_{\mathcal{R}}A$ is simple
by hypothesis (2) above.
Then $B = \End {}_{\mathcal{R}}A$ is an intermediate division ring between $C \subseteq A$, and $\mathcal{R} = \End A_B$ by Theorem~\ref{th-JBT}. 
 Since $\mathcal{I} \into
S$ and ${}_V\mathcal{I}$ is flat,
it follows from $A \o_V S \cong \End A_C$ that $\End A_B \cong A \o_V \mathcal{I}$ via the mapping above.
Note that $\mathcal{I} \subseteq \End A_B \cap S = \End {}_CA_B$
and let $Q$ be the cokernel. Since $A \o_V \mathcal{I} \cong \mathcal{R} \cong A \o_V \End {}_CA_B$
it follows  that $A \o_V Q = 0$. Since $A_V$ is faithfully flat,
$Q = 0$, whence 
$\mathcal{I} = \End {}_CA_B$. 
Finally, $\End A_B$ is isomorphic to an $A$-$C$-bimodule direct
summand of $A^N$, since ${}_V\mathcal{I} \oplus * \cong V^N$ for some
$N$, to which we apply the functor
${}_AA_C \o_V -$.  Since $A_B$ is finite free, it follows from Theorem~\ref{th-char}
that $A \supseteq B \supseteq C$ is left D3.  
\end{proof} 

If $A$ or $V$ is a division ring, the faithful flatness hypothesis
in the theorem is clearly satisfied.  
 In connection with this theorem we note the following criterion
for a depth three tower of division algebras.

\begin{prop}
Suppose $C \subseteq B \subseteq A$ is a tower 
of division rings where the right vector
space $A_B$ has basis $\{ a_1, \ldots,
a_n \}$ such that
\begin{equation}
Ca_i \subseteq a_iB \ \ \ (i = 1, \ldots,n) 
\end{equation}
Then $A \| B \| C$ is left D3.
\end{prop}
\begin{proof}
It is easy to compute that 
$x \o_B 1 = \sum_i a_i \o_B a_i^{-1} \beta_i(x)$ for all $x \in A$. Here  $\beta_i$ is the rank one projection
onto the right $B$-span of the basis element $a_i$
along the span of $a_1, \ldots, \hat{a_i},
\ldots, a_n$,
and $a_i^{-1} \o_B a_i \in (A \o_B A)^C$
for each $i$.  Of course, $\beta_i \in
\End {}_CA_B$, so $A \| B$ is left
D3.  
\end{proof}

We may similarly prove that the tower is
rD3 if ${}_BA$ has basis $\{ a_i \}$ satisfying
$a_i C \subseteq B a_i$. 
When $B = C$ we deduce the following criterion for a depth two
subalgebra pair of division rings. For example,
 the real quaternions $A = \H$,
 and subring  $B = \C$ meet this criterion.
\begin{cor}
Suppose $B \subseteq A$ is a subring pair
of division rings where the left vector
space ${}_BA$ has basis $\{ a_1, \ldots, a_n \}$ such that
\begin{equation}
a_i B = Ba_i \ \ \ (i = 1, \ldots,n).
\end{equation}
Then $A \| B$ is depth two.
\end{cor} 

We remark that 
if the centralizer $V$ of a depth two proper
extension $A \| C$ is contained in $C$ (as in the example $C = \C$ and $A = \H$ just mentioned above), then $\End {}_CA_C$ is a skew Hopf algebra
over the commutative base ring $V$ \cite{LK2007b}.  Any
intermediate ring $B$ of $A \| C$, for which $A \| B$ is D2,
has skew Hopf algebra $\End {}_BA_B$ over $R = A^B$ for the same reason, since
$R \subseteq V$ $ \subseteq $ $C \subseteq B$. It is interesting
to determine under what conditions these are skew Hopf subalgebras,
i.e., the antipodes are compatible under the sub-$R$-bialgebroid
structures.    

%%%%%%%%%%%%%%%%%%%%%%%%%%%%%%%%%%%%%%%%%%
\section{Application to field theory}

Given a separable finite field extension $F \subseteq E$
Szlach\'anyi shows that there is a Galois connection between
intermediate fields and weak Hopf subalgebras of $\End E_F$.
A \textit{weak Hopf algebra} $H$ the reader will recall from the already classic \cite{BNS}
is a weakening of the notion of Hopf algebra to include  certain
non-unital coproducts, non-homomorphic counits with
weakened antipode equations. There are certain canonical
coideal subalgebras $H^L$ and $H^R$ that are separable algebras
and anti-isomorphic copies of one another via the antipode.  
Nikshych and Etingof \cite{EN}
have shown that $H$ is a Hopf algebroid over the separable algebra
$H^L$, and conversely the author and Szlach\'anyi \cite{KS} have shown 
that Hopf algebroids over a separable algebra are weak Hopf algebras.  
Let's revisit one of the important, motivating examples. 

\begin{example}
\begin{rm}
Let $\mathcal{G}$ be a finite groupoid
with $x, y \in \mathcal{G}_{\rm obj}$ the objects
and $g,h \in \mathcal{G}_{\rm arrows}$ the invertible
arrows (with sample elements). Let $s(g)$ and $t(g)$
denote the source and target objects of the arrow $g$.  
 Suppose $k$ is a field.
Then the groupoid algebra $H = k\mathcal{G}$ (defined like
a quiver algebra, where $gh = 0$ if  $t(h) \neq s(g)$)
is a weak Hopf algebra with coproduct $\cop(g) = g \o_k g$,
counit $\eps(g) = 1$, and antipode $S(g) = g^{-1}$.   Since the identity is $1_H = \sum_{x \in \mathcal{G}_{\rm obj}} \id_x$,
we see that $\cop(1_H) \neq 1_H \o 1_H$
if $\mathcal{G}_{\rm obj}$ has two or more
objects. Notice too that $\eps(gh) \neq \eps(g) \eps(h)$
if $gh = 0$. 

The Hopf algebroid structure has total algebra $H$, and has base algebra
the separable
algebra $k\mathcal{G}_{\rm obj}$, which is a product algebra $k^N$
where $N = |\mathcal{G}_{\rm obj}|$. The source and target maps
of the Hopf algebras $s_L, t_L: R \rightarrow H$ are simply
$s_L = t_L: x \mapsto \id_x$.  The resulting bimodule
structure ${}_RH_R = {}_{s_L, \, t_L}H$ is given by $x \cdot g \cdot y = 
g$
if $x = y = t(g)$, $0$ otherwise. The coproduct is $\cop(g) = g\o_R g$,
counit $\eps(g) = t(g)$,  and antipode $S(g) = g^{-1}$.  
This defines a Hopf algebroid in the sense of Lu and Xu.
That this is also a Hopf algebroid in the sense of B\"ohm-Szlach\'anyi
may be seen by defining a right bialgebroid structure on $H$
via the counit $\eps_r(g) = s(g)$.   

If $\mathcal{G}$ is the finite set $\{ \underline{1}, \ldots, \underline{n} \}$ with singleton hom-groups, suggestively denoted by $\Hom (\underline{i},
\underline{j}) = \{ e_{ji} \}$ for all $1 \leq i,j \leq n$,
the groupoid algebra considered above is the full matrix algebra $H \cong M_n(k)$
and $R$ is subalgebra of diagonal matrices.   Note that the projection $\Pi^L$ ($ = \eps_t$ in \cite{EN}) defined as $\Pi^L(x) = \eps(1\1 x)1\2$ is given here by
$e_{ij} \mapsto e_{ii}$. Similarly,
$\Pi^R(e_{ij}) = e_{jj}$.  
\end{rm}
\end{example}

In \cite{Sz}, Szlach\'anyi shows that although Hopf-Galois separable field extensions
do not have a universal Hopf algebra as ``Galois quantum group,''
they have a universal weak Hopf algebra
or ''Galois quantum groupoid.''  For example,
the field $E = \Q(\sqrt[4]{2})$ is a four dimensional separable
extension of $F = \Q$ which is Hopf-Galois with respect to two
non-isomorphic Hopf algebras, $H_1$ and $H_2$ \cite{GP}.  However,
the endomorphism ring $\End E_F$ is then a smash product in two
ways, $E \# H_i$, $i = 1,2$, and is a weak Hopf algebra over the
separable $F$-algebra $E$.  It is universal in a category of weak
Hopf algebras viewed as left bialgebroids \cite[Theorem 2.2]{Sz}, with modifications to
the definition of the arrows resulting (see \cite[Prop.\ 1.4]{Sz}
for the definition of weak left morphisms of weak bialgebras).
The separable field extensions that are Hopf-Galois may then
be viewed as being weak Hopf-Galois with a uniqueness property.

The following corollary addresses an unanswered question in
\cite[Section 3.3]{Sz}. Namely, there is a Galois connection
between intermediate fields $K \subseteq F \subseteq E$
of a separable (finite) field extension $E \| K$
and weak Hopf subalgebras of the weak Hopf algebra $\mathcal{A} := \End E_K$
that include $E$ as left multiplications.  The
correspondences are denoted by $$\mbox{\rm Sub}_{WHA/K}(\mathcal{A})
\stackrel{\rm Fix}{\longrightarrow} \mbox{\rm Sub}_{Alg/K}(E)$$
which associates to a weak Hopf subalgebra $W$ of $\End E_K$ the
subfield $$\mbox{\rm Fix}(W) = \{ x \in E | \forall \alpha \in W,
\alpha(x) = \alpha(1) x \},$$ in other words, $E^W$,
and the correspondence $$\mbox{\rm Sub}_{Alg/K}(E) \stackrel{\rm Gal}{\longrightarrow} \mbox{\rm Sub}_{WHA/K}(\mathcal{A})$$
 where the intermediate subfield $K \subseteq F \subseteq E$
gets associated to its Galois algebra $$\mbox{\rm Gal}(F) = \{ \alpha \in \mathcal{A} |
\forall x \in E, y \in F, \alpha(xy) = \alpha(x)y \}. $$
Clearly Gal($F$) = $\End E_F$.  
 
Szlach\'anyi \cite[3.3]{Sz} notes that Gal is a surjective
correspondence, since $F$ = Fix(Gal($F$) for each intermediate
subfield (e.g. since 
$E_F$ is a generator module, it is  balanced
by Morita's lemma).   Gal is indeed a one-to-one correspondence
by 
\begin{cor}
Gal and Fix are inverse correspondences between intermediate fields
of a separable field extension $E \| K$ and weak Hopf subalgebras of
the full linear endomorphism algebra $\End E_K$.
\end{cor}
\begin{proof}
We just need to apply the Jacobson-Bourbaki correspondence with a 
change of notation.  Before changing notation, first note
that if $A \supseteq B$ is a depth two extension where $B$ is
a commutative subring of the center of $A$, then the centralizer $A^B = A$
and the left bialgebroid $\End {}_BA_B = \End A_B$ over $A$. Indeed, a faithfully flat $B$-algebra
$A$ is depth two iff it is finite projective. If $A$ and $B$ are fields,
this reduces to: depth two extension $A \| B$ 
$\Leftrightarrow$ finite extension $A \| B$.   
If $A \| B$ is a Frobenius
extension (as are separable extensions of
fields), there is an antipode
on $\End {}_BA_B$ defined
in terms of the Frobenius homomorphism
(such as the trace map of a separable field
extension \cite{SL}) 
 and its dual bases \cite{BS}.
Now, changing notation, we have a bialgebroid $\End E_K$ over the
separable $F$-algebra $E$, or equivalently a weak  bialgebra --- which becomes a weak
Hopf algebra via an involutive 
 antipode given
in terms of the trace map and its 
dual bases \cite[eq.~(3.5)]{Sz}). 

Given a weak Hopf subalgebra $W$ of $\End E_K$ containing $\lambda(E)$,
it is automatically finite dimensional over $E$ and ${}_WE$ is simple
since a submodule is a $W$-stable ideal, but $E$ is a field.  Hence,
$W$ is a Galois subring and the Theorem~\ref{th-JBT} shows that $\End {}_WE \cong
E^W$ is an intermediate field $F$ between $K \subseteq E$,
such that $\End E_F = W$.  But Gal($F$) = $\End E_F$ has
been noted above.  Hence, Gal(Fix($W$) = $W$. 
\end{proof}

The only reason we need restrict ourselves
to \textit{separable} field extensions
above is to acquire a fixed base algebra
that is a separable algebra, so that
we acquire antipodes from Frobenius extensions, and Hopf
algebroids become weak Hopf algebras.
Let us be clear on what happens when we drop this hypothesis.  For the
purpose of the next corollary , we define a  sub-$R$-bialgebroid of bialgebroid $(H,R,s_L,t_L \cop, \eps)$ to be a subalgebra $V$ of the total algebra $H$ with
the same base algebra $R$, source $s_L$ and target $t_L$
maps having image within $V$, and $V$ is
a sub-$R$-coring of $(H,\cop, \eps)$.

\begin{cor}
Let $E \supseteq K$ be a finite field
extension.  Then the poset of intermediate subfields
is in Galois correspondence with the poset of 
sub-$E$-bialgebroids of $\End E_K$.
\end{cor}
\begin{proof}
This follows from the Jacobson-Bourbaki
correspondence, where intermediate
field $F \mapsto \End E_F$
with inverse, Galois subring $R \mapsto \End {}_RE$, with the same proof as in the
previous corollary. Note from the proof
of Jacobson-Bourbaki in the field context  that  
any subring of $\End E_K$ containing
$\lambda(E)$ is indeed of the form
$\End E_F$ for some intermediate $K \subseteq F \subseteq E$, and therefore
the left bialgebroid of the depth two
(= finite) field extension $F \subseteq E$, and sub-$E$-bialgebroid of $\End E_K$.   
\end{proof}  
  
The Jacobson-Bourbaki correspondence also exists between subfields
of a finite dimensional simple algebra $A$ and subalgebras of the linear
endomorphism algebra which contain left and right multiplications
\cite[sect.\ 12.3]{P}, a theorem related to the topic of Brauer group of a field.  By the same reasoning, we arrive at Galois correspondences
between subfields and bialgebroids over $A$.  Namely, let
$A^e$ denote the image of $A \o_F A^{\rm op}$ in the linear
endomorphism algebra $\End A_F$ via left and right multiplication
$x \o y \mapsto \lambda_x \circ \rho_y$, and $Z(A)$ denote
the center of $A$, which is  a field since $Z(A) \cong \End {}_{A^e}A$.  We
note that $\End A_E$ is   a bialgebroid
over $A$ for any intermediate
field $F \subseteq E \subseteq Z(A)$ with
 Lu structure \cite{Lu}, and a Hopf algebroid in the special case $E = Z(A)$
where $A$ becomes Azumaya so $A \o_E A^{\rm op} \cong \End A_E$. 
The proof is quite the same as above and therefore omitted. 

\begin{cor}
Let $A$ be a simple finite dimensional $F$-algebra.  Then the fields
that are intermediate to $F \subseteq Z(A)$ are in Galois correspondence
to the sub-$A$-bialgebroids of $\End A_F$.  In case
$A$ is a separable $F$-algebra, the intermediate fields
are in Galois correspondence to weak Hopf subalgebras of $\End A_F$.  
\end{cor}

%%%%%%%%%%%%%%%%%%%%%%%%%%%%%%%%%%

\section{An embedding theorem for finite depth Frobenius extensions}

In this section we define finite depth Frobenius extension using the notion of depth three tower  by choosing a suitable three-ring sub-tower of the Jones tower.
We show that this definition is consistent with previous definitions of finite depth for subfactors and free Frobenius extensions.  We show in Theorem~\ref{th-embed} that any finite depth Frobenius extension extends to a depth two
extension somewhere further along in its Jones tower.  

Suppose $M_{-1} \into M_0$ is a Frobenius extension; e.g.
a (type $II_1$) subfactor of finite index \cite{K,KN, KS}
or a Frobenius algebra in the tensor category of bimodules over $M_{-1}$.  Let $M_1$
denote its basic construction $M_1 = Me_1M$ which is isomorphic
to $\End M_N$ and to $M \otimes_N M$ where $M:= M_0$ and $N := M_{-1}$.
The ring extension $M_1 \| M$ is itself a Frobenius extension
with $M$-bimodule Frobenius homorphism $E_1: M_1 \to M$ defined by 
$E_1(e_1) = 1$.    
The Jones element $e_1$ maps isomorphically into the Frobenius homomorphism and
into the cyclic generator $1_M \o_N 1_M$.  Iterate this to obtain
the Jones tower,
\begin{equation}
\label{eq: vfr}
N = M_{-1} \into M = M_0 \into M_1 \into M_2 \into \cdots \into M_{n} \into \cdots, 
\end{equation}
e.g., $M_2 = M_1 e_2 M_1$ where $e_2$ maps into the Frobenius homomorphism
$E_1: M_1 \to M_0$ and into $1_{M_1} \otimes_{M} 1_{M_1}$, $1_{M_1} = 
\sum_i x_i e_1 y_i$.  Note that $M_n \cong M \o_N \cdots \o_N M$
($n + 1$ times $M$).  Each ring extension $M_n \| M_{n-1}$ is a Frobenius
extension by the endomorphism ring theorem (iterated). Each
composite ring extension $M_n \| M_{n-k}$ is a Frobenius extension
by composing Frobenius homorphisms, and $M_{n+k}$ is isomorphic
to the basic construction of $M_n \| M_{n-k}$ by  \cite[appendix]{KN2}. 
While the $e_i$ may  not be idempotents or projections, they satisfy
$e_i e_{i \pm 1} e_i = e_i$, $e_i y e_i = E_{i-1}(y) e_i = e_i E_{i-1}(y)$  and $e_ix = e_i E_i(e_i x)$ for all $y \in M_{i-1}$, 
$x \in M_i$ with Frobenius homomorphisms $E_i: M_i \to M_{i-1}$.  For more details
on the Jones tower over a Frobenius extension, please see \cite[section 6]{KS}
and \cite[chapter 3]{K}. 

\begin{definition}[Finite depth Frobenius extensions]
The Frobenius extension $N \into M$ is said to be of depth $n > 1$
if the composite tower $M_{n-2} \| M_{n-3} \| M_{-1}$ is a right or left depth three tower.
\end{definition}
The definition allows for the possibility of a depth $n$ extension
being at the same time depth $n+1$, something we note to be true below. In subfactor theory,
one speaks of depth $n$ subfactor as the least $n$ for which the relative commutant $M_n^N$
is a basic construction of the two previous
semisimple algebras in the derived tower,
which we introduce next. 

Let $M_i^N$ denote the centralizer of $N$ in $M_i$.  We note the derived
tower of eq.~(\ref{eq: vfr}) above:
\begin{equation}
N^N \into M^N \into {M_1}^N \into \cdots \into {M_{n-1}}^N \into {M_n}^N \into \cdots
\end{equation}
In classical subfactor theory, depth $n$ is characterized by the least
$n$ for which $M_n^N$ is isomorphic to the basic construction of $M_{n-1}^N$
over $M_{n-2}^N$.  We compute that this is so with our new definition,
which we also show to be consistent with the definition in \cite[3.1]{KN}
of depth $n$ free Frobenius extension.  

\begin{prop}
A depth $n$ Frobenius extension $N \into M$ has $n$-step
centralizer $M_n^N \cong M_{n-1}^N \otimes_{M_{n-2}^N} M_{n-1}^N$.
The Frobenius homomorphism $E_{n-1}$ has dual bases elements in
$M_{n-1}^N$.   
\end{prop}

\begin{proof}
Let $A = M_{n-2}$, $B = M_{n-3}$ and $C = M_{-1}$.  Then subtower
$A \| B \| C$ of the Jones tower~(\ref{eq: vfr}) is D3.  Whence
the $A$-$C$-bimodules $A \o_B A$ and $A$ are $H$-equivalent,
and their endomorphism rings are Morita equivalent.

Note that $$\End {}_A A \otimes_B A_C \cong \End {}_A( M_{n-1})_N \cong
M_n^N,$$
since $M_{n-1}$ is isomorphic to the basic construction of $M_{n-2} \| M_{n-3}$ and anti-isomorphic to $\End {}_BA$.
We use as well that the $C$-centralizer $(\End {}_BA)^C \cong \End {}_BA_C$. 
  
On the other hand, we note
$$\End {}_AA_C \cong A^C \cong M_{n-2}^N,$$
so we conclude $M_n^N$ and $M_{n-2}^N$ are Morita equivalent rings.
The Morita context bimodules are the $A$-$C$-bimodule hom-groups
$\Hom (A \o_B A, A)$ and $\Hom (A, A \o_B A)$.  In section~2 we saw
that first of all, 
$$ \Hom (A \o_B A, A) \cong \End {}_BA_C \cong M_{n-1}^N,$$
since $M_{n-1}$ is anti-isomorphic to the left endomorphism ring
of $M_{n-2} \| M_{n-3}$.  Since the Frobenius extension $A \| B$
satisfies $\End {}_BA \cong A \o_B A$, we also obtain that the $C$-centralizers, 
$$\End {}_BA_C \cong (A \o_B A)^C \cong \Hom (A, A \o_B A)$$
the last step following from section~2. In other words,
$M_{n-1}^N$ doubles as both of the Morita context bimodules between
$M_n^N$ and $M_{n-2}^N$.   

By Morita theory, it follows that the Morita context bimodules satisfy  
$ M_n^N \cong M_{n-1}^N \o_{M_{n-2}^N} M_{n-1}^N$. 

The last statement is proven as  \cite[Theorem 2.1, part (5)]{KK}
and \cite[Prop.\ 6.4]{KS}.  Toward this end, we note that $M_{n-1} \cong Ae_{n-1}A$ and $M_{n-1}^N \cong (A \o_B A)^C \cong \End {}_CA_B$.  Suppose a left D3 quasibases  given by $t_i = t^1_i e_{n-1} t^2_i$ and $\beta_i \in \End {}_CA_B = (\End A_B)^C$; 
with $\sum_j x_j \o_B y_j \in (A \o_B A)^A$  dual bases for
the Frobenius homomorphism $ E_{n-2}: A \to B$.  Then $E_{n-1}: M_{n-1}=
A e_{n-1}A \to A$ has dual bases  $t_i = t_i^1 e_{n-1} t_i^2$
and  $\sum_j \beta_i(x_j)e_{n-1} y_j$, both in $M_{n-1}^N$.
\end{proof}

We note that depth $n$ Frobenius extensions are depth $n+1$
as follows.

\begin{lemma}[Endomorphism ring lemma for D3 towers]
Suppose $A \| B \| C$ is D3 where $A \| B$ is a Frobenius extension.
Let $D = \End A_B$ and $D \| A$ denote the left regular mapping
$a \mapsto \lambda_a$.  Then the tower $D \| A \| C$ is D3.
\end{lemma}
\begin{proof}
Since $A \| B$ is a Frobenius extension, we note that
${}_DD_A \cong {}_D A \o_B A_A$ given by
$d \mapsto \sum_i d(x_i) \o_B y_i$ in the Frobenius system notation
 above.  Hence ${}_DD \o_A D_C \cong {}_DA \o_B A \o_B A_C$.
To ${}_AA \o_B A_C \oplus * \cong {}_AA^m_C$
 we apply the functor ${}_DA \o_B - $.  We obtain after substitution,
${}_DD \o_A D_C \oplus * \cong {}_DD^m_C$, i.e. the tower $D \| A \| C$
is rD3.  We have noted above eq.~(\ref{eq: left-right}) that
$D \| A \| C$ is necessarily also $\ell$D3
since $D \| A$ is a Frobenius extension
\end{proof}

The next lemma argues as in Theorem~\ref{th-conv}.

\begin{lemma}[Tunneling Lemma]
Suppose the tower $A \| B \| C$ is rD3.  Let $C \rightarrow B$
factor through a Frobenius extension $D \rightarrow B$ such
that $A \cong \End B_D$.
  Then $A \| D \| C$ is rD3.
\end{lemma}
\begin{proof}
We have $A \cong B \o_D B$, from which we obtain $A \o_B A \o_B A \cong
A \o_D A$.  Tensoring $A \o_B A \oplus * \cong A^n$ by $A \o_B -$,
we obtain this fact.
\end{proof}

The next theorem shows that a depth $n = 2^m + 1 = 3, 5, 9, 17, \ldots$ Frobenius extension
$N \into M$ is embedded in (``is a factor of'') the depth two extension
$N \into M_{n-2}$.
  
\begin{theorem}
\label{th-embed}
If $N \into M$ is a depth $n$ Frobenius extension,
where $n = 2^m + 1$ for a positive integer $m$, then
$N \into M_{n-2}$ is a depth two Frobenius extension.
\end{theorem}

\begin{proof}
By hypothesis, the Jones subtower $M_{n-2} \| M_{n-3} \| M_{-1}$ is D3.
Since $M_{n-2} \cong M_{n-3} \o_{M_{n-4}} M_{n-3}$, it follows from
the tunneling lemma that $M_{n-2} \| M_{n-4} \| M_{-1}$ is D3.  Iterating use
of the lemma, also $M_{n-2} \| M_{n-6} \| M_{-1}$ is D3.
Continuing $m-1$ steps, $M_{n-2} \| M_{2^{m-1} - 1} \| M_{-1}$ is
D3.  But $M_{n-2} \cong M_{2^{m-1} - 1} \o_N M_{2^{m-1} - 1}$
since $n - 2 = 2^m -1$ (and recalling $M_i = M^{\otimes_N \, i+1}$).  By Theorem~\ref{th-conv} $M_{n-2} \| N$ is D2.
\end{proof}

Suppose a Frobenius extension $N \into M$ is depth $n$,
where $2^{m-1} < n \leq 2^m$.  Then by the endomorphism ring lemma
this extension is also depth $2^m + 1$, so
$N \into M_{2^m -1}$ is D2 by the theorem. This proves: 

\begin{cor}
Any finite depth Frobenius extension $N \into M$ embeds
into a depth two extension of $N$ in  an $n$'th 
iterated endomorphism ring $M_n$.
\end{cor}

This result may be viewed as an algebraic version of
a result in \cite{NV} and an answer to a question in \cite[appendix]{KN2}.  
The theory above seems to indicate that one of a variety of  extensions
of Jacobson-Bourbaki correspondence to pairs of simple algebras
by the Japanese school of ring theory (Hirata, M\"uller, Onodera,  Sugano, 
Szeto, Tominaga, and others)
would adapt via depth two extensions and depth three towers to
 an algebraic version of the Galois theory
for subfactors in Nikshych and Vainerman
\cite{NV}. This will be investigated in
a future paper.   

Since group algebras and their finite index subgroups  form Frobenius extensions, we pose the following question in character theory and group theory in extension of the discussion and results in sections~3 and~8 of this paper.

Question:  what precisely are the group-theoretic conditions on a  subgroup of a finite group  $H < G$ that its 
Frobenius algebra extension $N = \C [H] \subseteq M = \C [G]$  be depth $n$?

\end{document}